\documentclass[11pt]{article}
\usepackage{cite}
\usepackage{mathrsfs}
\usepackage{amsfonts}
\usepackage{amsmath}
\usepackage{amsfonts,amssymb,color}
\usepackage{dsfont}
\usepackage{curves}
\usepackage{mathrsfs}
\usepackage{pifont}
\usepackage{amssymb}
\usepackage{latexsym,amsmath,amssymb,amsfonts,epsfig,graphicx,cite,psfrag}
\usepackage{eepic,color,colordvi,amscd}

\newtheorem{theorem}{Theorem}[section]

\newtheorem{lemma}{Lemma}[section]

\newtheorem{claim}{Claim}[section]
\newtheorem{conjecture}{Conjecture}[section]

\newcommand{\qed}{\hfill\rule{0.5em}{0.809em}}

\def\emptyset{\mbox{{\rm \O}}}
\def\bar{\overline}

\def\qed{\hfill \rule{4pt}{7pt}}

\def\pf{\noindent {\it Proof. }}

\begin{document}

 \title{Perfect divisibility and coloring of some fork-free graphs\thanks{Partially supported by NSFC projects 11931006 and 12126339}}
  \author{Di Wu$^{1,}$\footnote{Email: 1975335772@qq.com},  \;\; Baogang  Xu$^{1,}$\footnote{Email: baogxu@njnu.edu.cn OR baogxu@hotmail.com.}\\\\
 	\small $^1$Institute of Mathematics, School of Mathematical Sciences\\
 	\small Nanjing Normal University, 1 Wenyuan Road,  Nanjing, 210023,  China}
 \date{}

 \maketitle
\begin{abstract}
A $hole$ is an induced cycle of length at least four, and an odd hole is a hole of odd length. A {\em fork} is a graph obtained from $K_{1,3}$ by subdividing an edge once.  An {\em odd balloon} is a graph obtained from an odd hole by identifying respectively two consecutive vertices with two leaves of $K_{1, 3}$. A {\em gem} is a graph that consists of a $P_4$ plus a vertex adjacent to all vertices of the $P_4$. A {\em butterfly} is a graph obtained from two traingles by sharing exactly one vertex. A graph $G$ is perfectly divisible if for each induced subgraph $H$ of $G$, $V(H)$ can be partitioned into $A$ and $B$ such that $H[A]$ is perfect and $\omega(H[B])<\omega(H)$.  In this paper, we show that (odd balloon, fork)-free graphs are perfectly divisible (this generalizes some results of Karthick {\em et al}). As an application, we show that $\chi(G)\le\binom{\omega(G)+1}{2}$ if $G$ is (fork, gem)-free or (fork, butterfly)-free.

\begin{flushleft}
{\em Key words and phrases:} perfect divisibility, fork, chromatic number\\
{\em AMS 2000 Subject Classifications:}  05C15, 05C75\\
\end{flushleft}

\end{abstract}

\newpage

\section{Introduction}
All graphs considered in this paper are finite and simple. We follow \cite{BM08} for undefined notations and terminology.  We use $P_k$ and $C_k$ to denote a {\em path} and a {\em cycle} on $k$ vertices respectively. Let $G$ be a graph, $v\in V(G)$, and let $X$ and $Y$ be two subsets of $V(G)$. We say that $v$ is {\em complete} to $X$ if $v$ is adjacent to all vertices of $X$, and say that $v$ is {\em anticomplete} to $X$ if $v$ is not adjacent to any vertex of $X$. We say that $X$ is complete (resp. anticomplete) to $Y$ if each vertex of $X$ is complete (resp. anticomplete) to $Y$. If $1<|X|<|V(G)|$ and every vertex of $V(G)\setminus X$ is either complete or anticomplete to $X$, then we call $X$ a {\em homogeneous set}. We use $G[X]$ to denote the subgraph of $G$ induced by $X$, and call $X$ a {\em clique} if $G[X]$ is a complete graph. The {\em clique number} $\omega(G)$ of $G$ is the maximum size taken over all cliques of $G$.

For $v\in V(G)$ and $X\subseteq V(G)$, let $N_G(v)$ be the set of vertices adjacent to $v$, $d_G(v)=|N_G(v)|$, $N_G[v]=N_G(v)\cup \{v\}$, $M_G(v)=V(G)\setminus N_G[v]$, and let $N_G(X)=\{u\in V(G)\setminus X\;|\; u$ has a neighbor in $X\}$ and $M_G(X)=V(G)\setminus (X\cup N_G(X))$. If it does not cause any confusion, we usually omit the subscript $G$ and simply write $N(v), d(v), N[v], M(v), N(X)$ and $M(X)$.  Let $\delta(G)$ denote the minimum degree of $G$.

Let $G$ and $H$ be two vertex disjoint graphs. The {\em union} $G\cup H$ is the graph with $V(G\cup H)=V(G)\cup (H)$ and $E(G\cup H)=E(G)\cup E(H)$. The {\em join} $G+H$ is the graph with $V(G+H)=V(G)+V(H)$ and $E(G+H)=E(G)\cup E(H)\cup\{xy\;|\; x\in V(G), y\in V(H)$$\}$. The union of $k$ copies of the same graph $G$ will be denoted by $kG$.  We say that $G$ induces $H$ if $G$ has an induced subgraph isomorphic to $H$, and say that $G$ is $H$-free otherwise. Analogously, for a family $\cal H$ of graphs, we say that $G$ is ${\cal H}$-free if $G$ induces no member of ${\cal H}$.

For $u, v\in V(G)$, we simply write $u\sim v$ if $uv\in E(G)$, and write $u\not\sim v$ if $uv\not\in E(G)$.  A {\em hole} of $G$ is an induced cycle of length at least 4, and a {\em $k$-hole} is a hole of length $k$. A $k$-hole is called an {\em odd hole} if $k$ is odd, and is called an {\em even hole} otherwise. An {\em antihole} is the complement of some hole. An odd (resp. even) antihole is defined analogously.

Let $k$ be a positive integer, and let $[k]=\{1, 2, \ldots, k\}$. A $k$-{\em coloring} of $G$ is a mapping $c: V(G)\mapsto [k]$ such that $c(u)\neq c(v)$ whenever $u\sim v$ in $G$. The {\em chromatic number} $\chi(G)$ of $G$ is the minimum integer $k$ such that $G$ admits a $k$-coloring.  It is certain that $\chi(G)\ge \omega(G)$. A {\em perfect graph} is one such that  $\chi(H)=\omega(H)$ for all of its induced subgraphs $H$. The famous {\em Strong Perfect Graph Theorem}\cite{CRSR06} states that a graph is perfect if and only if it induces neither an odd hole nor an odd antihole. A graph is {\em perfectly divisible} \cite{HCT18} if for each of its induced subgraph $H$, $V(H)$ can be partitioned into $A$ and $B$ such that $H[A]$ is perfect and $\omega(H[B])<\omega(H)$. By a simple induction on $\omega(G)$ we have that $\chi(G)\le \binom{\omega(G)+1}{2}$ for each perfectly divisible graph.

\begin{figure}[htbp]\label{fig-1}
	\begin{center}
		\includegraphics[width=12cm]{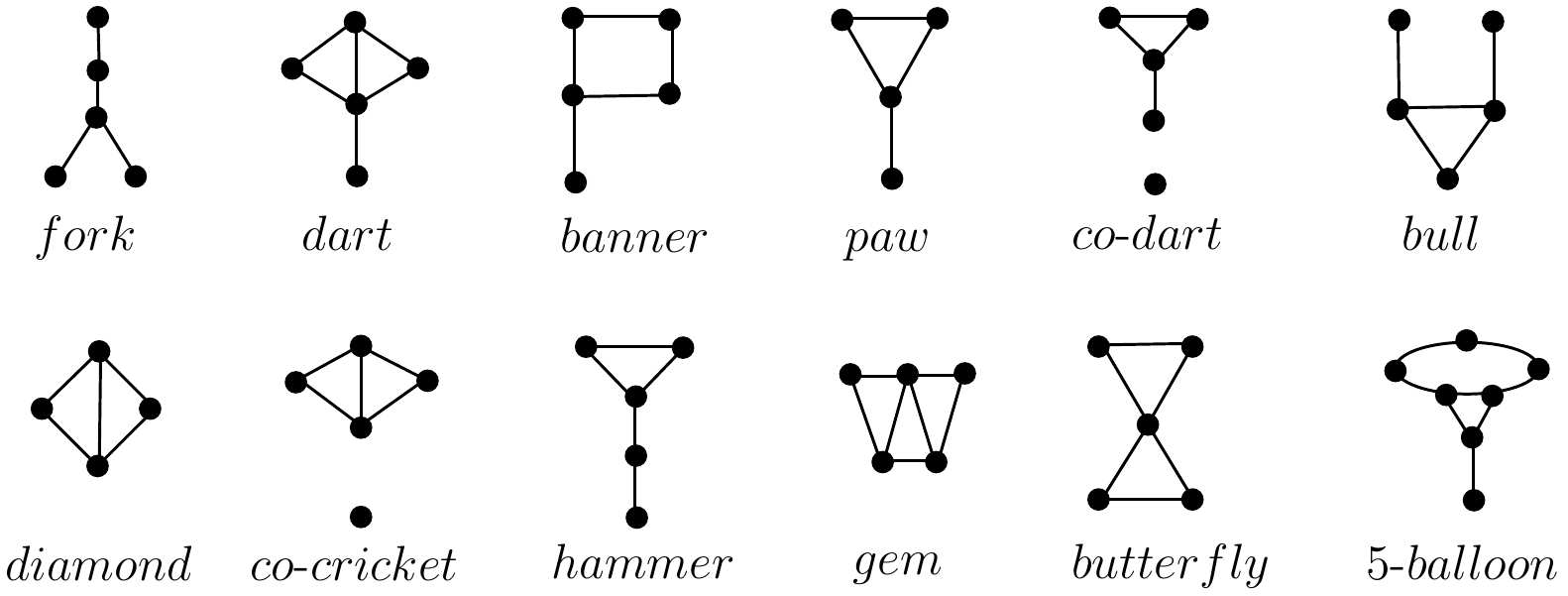}
	\end{center}
	\vskip -15pt
	\caption{Illustration of fork and some forbidden configurations.}
\end{figure}

A {\em fork} is the graph obtained from $K_{1,3}$ (usually called {\em claw}) by subdividing an edge once (see Figure \ref{fig-1}). The class of claw-free graphs is a subclass of fork-free graphs. It is natural to see what properties of claw-free graphs are also enjoyed by fork-free graphs. Sivaraman  proposed the following still open conjecture (see \cite{KKS22}).

\begin{conjecture}\label{conj-S}
   The class of fork-free graphs is perfectly divisible.
\end{conjecture}

Here we list some results related to (fork,$H$)-free graphs for some small graph $H$, and refer the readers to \cite{KKS22,CHKK21} for more information on Conjecture \ref{conj-S} and related problems.

A {\em dart} is the graph $K_1+(K_1\cup P_3)$, a {\em banner} is a graph consisting of a $C_4$ with one pendant edge, a {\em paw} is a graph obtained from $K_{1, 3}$ by adding an edge joining two of its leaves, a {\em co-dart} is the union of $K_1$ and a paw, a {\em bull} is a graph consisting of a triangle with two disjoint pendant edges, a {\em diamond} is the graph $K_1+P_3$, a {\em co-cricket} is the union of  $K_1$ and a diamond, a {\em hammer} is a graph obtained
by identifying one vertex of a triangle and one end vertex of a $P_3$, a {\em gem} is the graph $K_1+P_4$, a {\em butterfly} is the graph $K_1+2K_2$. A {\em balloon} is a graph obtained from a hole by identifying respectively two consecutive vertices with two leaves of $K_{1, 3}$. An $i$-{\em balloon}  is a balloon such that its hole has $i$ vertices. An $i$-balloon is called an {\em odd balloon}  if $i$ is odd. (See Figure \ref{fig-1} for these configurations)

Let ${\cal G}$ be a family of graphs. If there exists a function $f$ such that $\chi(G)\le f(\omega(G))$ for each graph $G\in {\cal G}$, then we say that ${\cal G}$ is $\chi$-{\em bounded} with binding function $f$ \cite{gyarfas1}.

In \cite{CHKK21, CCS20, RB93}, the authors proved that (fork, $H$)-free graphs are linearly $\chi$-bounded when $H\in\{C_4, K_4,$ diamond, $K_3\cup K_1$, paw, $\bar{P_3\cup 2K_1}, K_5-e$, antifork$\}$, where an {\em antifork} is the complement graph of a fork.

Kierstead and Penrice \cite{KP94} showed that for any tree $T$ of radius two, the class of $T$-free graphs is $\chi$-bounded. As a corollary, we see that fork-free graphs are $\chi$-bounded. In\cite{KKS22, RS04}, the authors asked that whether there exists a polynomial (in particular, quadratic) $\chi$-binding function for fork-free graphs. Very recently, Liu {\em et al}  \cite{LSWY21} answered the problem and proved that $\chi(G)\le7.5\omega^2(G)$ for every fork-free graph $G$. This is a significant progress on Conjecture~\ref{conj-S} that motivates us to pay more attention to the perfect divisibility of fork-free graphs.

Karthick, Kaufmann and Sivaraman \cite{KKS22} proved that Conjecture~\ref{conj-S} is true on (fork, $F$)-free graphs when $F\in$ \{$P_6$, co-dart, bull\}. They also showed that the classes of (fork, $H$)-free graphs, when $F\in\{$ dart, banner, co-cricket$\}$, are either claw-free or perfectly divisible.

In this paper, we generalize some results of Karthick {\em  et al} in \cite{KKS22} and prove the following theorem.

\begin{theorem}\label{fork}
  $($fork, odd balloon$)$-free graphs are perfectly divisible.
\end{theorem}

As a corollary, the class of (fork, $F$)-free graphs is perfectly divisible when $F\in \{P_6$, co-dart, bull, hammer, odd hole\}. As an application of Theorem~\ref{fork}, we also show that

\begin{theorem}\label{gem}
	$\chi(G)\le\binom{\omega(G)+1}{2}$ if $G$ is $($fork, gem$)$-free.
\end{theorem}

\begin{theorem}\label{butterfly}
	$\chi(G)\le\binom{\omega(G)+1}{2}$ if $G$ is $($fork, butterfly$)$-free.
\end{theorem}

We will prove Theorem~\ref{fork} in Section 2, prove Theorem~\ref{gem} in Section 3, and  prove Theorem~\ref{butterfly} in Section 4.

\section{(fork, odd balloon)-free graphs }

We prove Theorem~\ref{fork} in this section.
Let $X$ be a subset of $V(G)$. We say that $X$ is {\em anticonnected} in $G$ if the subgraph induced by $X$ in the complement $\bar{G}$ of $G$ is connected. An {\em antipath} in $G$ is the complement of a path of $\bar{G}$. A class of graphs ${\cal C}$ is {\em hereditary} if $G\in {\cal C}$ then all induced subgraphs of $G$ are also in ${\cal C}$. The following two lemmas will be used in our proof.

\begin{lemma}\label{pb-2}{\em \cite{CS19}}
A minimal nonperfectly divisible graph does not have homogeneous sets.
\end{lemma}

\begin{lemma}\label{pb-1}{\em \cite{KKS22}}
Let ${\cal C}$ be a hereditary class of graphs. Suppose that every graph $H\in{\cal C}$ has a vertex $v$ such that $H[M_H(v)]$ is perfect. Then every graph in ${\cal C}$ is perfectly divisible.
\end{lemma}

With a similar technique used by Karthick {\em et al} in \cite{KKS22}, we prove the next lemma.

\begin{lemma}\label{2-3}
Let $G$ be a connected fork-free graph. If $G$ is minimal nonperfectly divisible, then for each vertex $v\in V(G)$, $M(v)$ contains no odd antihole except $C_5$.
\end{lemma}

\pf Suppose to its contrary, let $G$ be a minimal nonperfectly divisible graph, and let $v\in V(G)$ such that $G[M(v)]$ contains an odd antihole induced by $X_0$ with $|X_0|\ge7$. Let $X_0=\{v_1, v_2,\ldots, v_n\}$ such that $v_i\sim v_j$ whenever $|i-j|\ne 1$ (indices are modulo
$n$).

We construct a sequence of vertex subsets $X_0, X_1,\dots, X_t,\cdots$ such that for each $i\ge 1$, $X_i$ is obtained from $X_{i-1}$ by adding one vertex $x_i$ that has a neighbor and a nonneighbor in $X_{i-1}$. Let $X$ be the maximal vertex set obtained in this way. By its maximality, we have that every vertex in $V(G)\setminus X$ is either complete or anticomplete to $X$. We will show that $X\ne V(G)$. Then $X$ is a homogeneous set which contradicts Lemma~\ref{pb-2}. Let $A=M(X_0)$. It is certain that $A\ne \emptyset$ as $v\in A$. Firstly, we will prove that

\begin{equation}\label{eqa-1-2}
\mbox{$N(A)$ is complete to $X_0$}.
\end{equation}
Suppose to the contrary that there exists a vertex $b\in N(A)$ that has a nonneighbor in $X_0$. Let $a$ be a neighbor of $b$ in $A$. Since $b\not\in A$, we have that $b$ has neighbors in $X_0$.

If $b$ has two consecutive neighbors in $X_0$, then there is some $i$ such that $v_i, v_{i+1}\in N(b)$ and $v_{i+2}\notin N(b)$, which implies a fork induced by $\{v_{i+2}, v_i, b, v_{i+1}, a\}$. This shows that $b$ does not have two consecutive neighbors in $X_0$.

Without loss of generality, suppose that $b\sim v_1$. Then $b\not\sim v_2$ and $b\not\sim v_n$.

If $b\sim v_3$, then $b\not\sim v_4$ implying that $b\sim v_5$ as otherwise $\{a,b,v_1,v_4,v_5\}$ induces a fork, and consequently $b\sim v_j$ for $j\in\{3,5,7,\dots,n-2\}$ and $b\not\sim v_k$ for $k\in\{2,4,6,\dots,n-1\}$. In particular, $b\sim v_{n-4}$ and $b\not\sim v_{n-1}$, which implies that $b\sim v_n$ as otherwise $\{a,b,v_{n-4},v_{n-1},v_n\}$ induces a fork, a contradiction. Therefore, $b\not\sim v_3$. Consequently, we have that $b\sim v_4$ to forbid a fork on $\{a,b,v_1,v_3,v_4\}$, and $b\not\sim v_5$ as $b$ has no two consecutive neighbors in $X_0$, and $b\sim v_6$ as otherwise $\{a,b,v_1,v_5,v_6\}$ induces a fork. But now, $\{a,b,v_6,v_2,v_3\}$ induces a fork, a contradiction. This proves (\ref{eqa-1-2}).

Next, we prove that

\begin{equation}\label{eqa-1-3}
\mbox{for $i\ge 0$, $X_i$ is complete to $N(A)$, and $X_i\cap (A\cup N(A))=\emptyset$}.
\end{equation}

By (\ref{eqa-1-2}), we have that (\ref{eqa-1-3}) holds for $i=0$. Suppose $i\ge 1$, and suppose that (\ref{eqa-1-3}) holds for $X_{i-1}$. Recall that $X_i=X_{i-1}\cup\{x_i\}$, where $x_i$ has a neighbor and a nonneighbor in $X_{i-1}$.  Suppose that $X_i$ is not complete to $N(A)$, and let $b\in N(A)$ be a nonneighbor of $x_i$.  Since $x_i$ has a neighbor in $X_{i-1}$ and $N(A)\cap X_{i-1}=\emptyset$, we have that $x_i\notin A$. Since $x_i$ has a nonneighbor in $X_{i-1}$ and $N(A)$ is complete to $X_{i-1}$, we have that $x_i\notin N(A)$. This shows that $X_i\cap (A\cup N(A))=\emptyset$ by induction.

Since $x_i\notin A$, $x_i$ has a neighbor in $X_0$, say $v_1$ by symmetry.
Notice that $X_i$ is anticonnected for each $i$. Let $P$ be a shortest antipath from $x_i$ to $v_1$. Then, $P$ has at least three vertices. Suppose that $P=w_1w_2\ldots w_t$, where $x_i=w_1$ and $v_1=w_t$. Since $X_{i-1}$ is complete to $N(A)$, and since $w_2, w_3\in X_{i-1}$, we have that $b\sim w_2$ and $b\sim w_3$, and so $\{x_i, w_2, b, w_3, a\}$ induces a fork which is a contradiction. Therefore, $x_i$ has no nonneighbor in $N(A)$, which means that $X_i$ is complete to $N(A)$.  This proves (\ref{eqa-1-3}).

It follows directly, from (\ref{eqa-1-3}), that $X\cap (A\cup N(A))=\emptyset$. Since $A\ne\emptyset$, we have that $X\ne V(G)$, which implies that $X$ is homogeneous set. This completes the proof of Lemma~\ref{2-3}.  \qed

\medskip

Let $G$ be a minimal nonperfectly divisible fork-free graph, and let $v$ be a vertex of $G$. By Lemma~\ref{pb-1}, $G[M(v)]$ is not perfect and must contain some odd hole by Lemma~\ref{2-3}.

\begin{lemma}\label{2-4}
	Let $G$ be a  fork-free graph, and $C=v_1v_2\dots v_nv_1$ an odd hole contained in $G$. If there exist two adjacent vertices $u$ and $v$ in $V(G)\setminus V(C)$ such that $u$ is not anticomplete to $V(C)$ but $v$ is, then $N(u)\cap V(C)=\{v_j, v_{j+1}\}$, for some $j\in \{1,2,\dots,n\}$, or $N(u)\cap V(C)= V(C)$.
\end{lemma}
\pf Let $C=v_1v_2\dots v_nv_1$ such that $n\ge5$ and $n$ is odd. By symmetry, we may assume that $u\sim v_1$. By the definition of $v$, we have that $u\sim v$ and $v$ is anticomplete to $V(C)$. Suppose that $N(u)\cap V(C)\not\in\{\{v_1,v_2\},\{v_1,v_n\}\}$.

If $u\not\sim v_2$ and $u\not\sim v_n$, then $\{v_1,v_2,v_n,u,v\}$ induces a fork, a contradiction.

If $u$ has exactly one neighbor in $\{v_2,v_n\}$, we suppose by symmetry that $u\sim v_n$, and let $t$ be the smallest integer in $\{3,4,\dots,n-1\}$ such that $u\sim v_t$, then $\{v_1,v_{t-1},v_t,u,v\}$ induces a fork when $t= n-1$, and $\{v_{t-1},v_t,v_n,u,v\}$ induces a fork when $t\ne n-1$ , both are contradictions.

Suppose that $u\sim v_2$ and $u\sim v_n$. If $u$ is not complete to $\{v_3,v_4,\dots v_{n-2}\}$, let $t'$ be the smallest integer in $\{3,4,\dots,n-2\}$ such that $u\not\sim v_{t'}$, then $\{v_{t'-1},v_{t'},v_n,u,v\}$ induces a fork, a contradiction. This shows that $u$ is complete to $\{v_3,v_4,\dots v_{n-2}\}$. Moreover, $u\sim v_{n-1}$ as otherwise $\{v_2,v_{n-1},v_n,u,v\}$ induces a fork, a contradiction. Therefore, $u$ is complete to $V(C)$. This proves Lemma~\ref{2-4}.  \qed

\begin{lemma}\label{2-5}
Let $G$ be a  connected minimal nonperfectly divisible fork-free graph, $v\in V(G)$, and $C$ an odd hole contained in $G[M(v)]$. Then $N(M(C))$ is not complete to $V(C)$.
\end{lemma}
\pf Suppose to its contrary, let $C=v_1v_2\dots v_nv_1$ such that $n\ge5$ and $n$ is odd,  $A=M(C)$, $X=N(A)$ that is complete to $V(C)$. Since $v\in A$, we have that $A\ne\emptyset$, and thus $X\ne\emptyset$ as $G$ is connected. Let $Z=N(C)\setminus X$. It is certain that $V(G)=V(C)\cup X\cup A\cup Z$.

Since $V(C)$ is anticomplete to $A$, and is complete to $X$ by our assumption, we have that $Z\ne\emptyset$ as otherwise $V(C)$ is a homogeneous set, contradicts Lemma \ref{pb-2}. Moreover, by the definition of $Z$, $Z$ is anticomplete to $A$.

If $Z$ is complete to $V(C)$, then $V(C)$ is a homogeneous set. If $Z$ is complete to $X$, then $V(C)\cup Z$ is a homogeneous set. Both contradict Lemma \ref{pb-2}. Therefore, we have that

\begin{equation}\label{2-4-1}
\mbox{$Z$ is neither complete to $V(C)$ nor complete to $X$.}
\end{equation}

Let $Z'$ be the set of vertices in $Z$ which are complete to $X$, let $Z''=Z\setminus Z'$. We will show that

\begin{equation}\label{2-4-2}
\mbox{$Z''\ne\emptyset$ and $Z''$ is complete to $V(C)$, and $Z'\ne\emptyset$.}
\end{equation}

By (\ref{2-4-1}), $Z$ is not complete to $X$, which implies that $Z''\ne\emptyset$. Since no vertex in $Z''$ is complete to $X$, each vertex $z$ of $Z''$ has a nonneighbor, say $x$, in $X$. Let $a\in A$ such that $a\sim x$, and suppose by symmetry that $v_1\sim z$. Then $z$ is complete to $V(C)\setminus \{v_2,v_n\}$ as otherwise $\{v',v_1,z,x,a\}$ induces a fork, where $v'$ is a nonneighbor of $z$ in $V(C)\setminus \{v_2,v_n\}$, a contradiction. Consequently, $z$ is complete to $\{v_2,v_n\}$ as otherwise $\{v_2,v_{n-1},x,z,a\}$ induces a fork when $z\not\sim v_2$, and $\{v_3,v_n,x,z,a\}$ induces a fork when $z\not\sim v_n$, both are contradictions. Therefore, $z$ is complete to $V(C)$, which implies that $Z''$ is complete to $V(C)$. If $Z'=\emptyset$, then $Z=Z''$ is complete to $V(C)$, contradicts (\ref{2-4-1}). Hence, $Z'\ne\emptyset$. This proves (\ref{2-4-2}).

Let $T$ be the set of vertices in $Z'$ which are complete to $Z''$, let $T'=Z'\setminus T$. Next, we prove that

\begin{equation}\label{2-4-3}
\mbox{$T\ne\emptyset$ and $T'\ne\emptyset$.}
\end{equation}

If $T'=\emptyset$. Then $Z'$ is complete to $Z''$, which implies that $V(C)\cup Z'$ is a homogeneous set as $V(C)\cup Z'$ is complete to $X\cup Z''$ and anticomplete to $A$, contradicts Lemma \ref{pb-2}. This shows that $T'\ne\emptyset$.

Suppose that $T=\emptyset$. Then each vertex in $Z'$ has a nonneighbor in $Z''$. By (\ref{2-4-1}) and (\ref{2-4-2}), we have that $Z'$ is not complete to $V(C)$. Let $z'\in Z'$ which has a nonneighbor in $V(C)$, say $v_1$, and let $z''\in Z''$ be a nonneighbor of $z'$. By the definition of $Z''$, $z''$ has a nonneighbor, say $x''$, in $X$. Let $a''$ be a neighbor of $x''$ in $A$. Since $X$ is complete to $V(C)\cup Z'$, we have that $x''\sim v_1$ and $x''\sim z'$ . By (\ref{2-4-2}), $Z''$ is complete to $V(C)$, in particular, $z''\sim v_1$. Therefore, $\{z'',v_1,z',x'',a''\}$ induces a fork, a contradiction. This proves (\ref{2-4-3}).

Suppose that $T$ is not complete to $T'$. Let $t\in T$ and $t'\in T'$ such that $t\not\sim t'$. Then $t'$ has a nonneighbor $b$ in $Z''$, $b$ has a nonneighbor $b'$ in $X$, $b'$ has a neighbor $c$ in $A$. Since $X$ is complete to $Z'$ and $T$ is complete to $Z''$, we have that $b'\sim t,b'\sim t'$ and $t\sim b$. Now, $\{t,t',b,b',c\}$ induces a fork, a contradiction. So, $T$ is complete to $T'$.

Finally, suppose that $T'$ is not complete to $V(C)$. By symmetry, we may assume that $d\in T'$ such that $d\not\sim v_1$. Notice that each vertex in $T'$ has a nonneighbor in $Z''$ and each vertex in $Z''$ has a nonneighbor in $X$. We suppose that $d'\in Z''$ such that $d'\not\sim d$. Consequently, we suppose that $d''\in X$ such that $d''\not\sim d'$. Recall $Z'$ is complete to $X$. It follows that $d\sim d''$. By (\ref{2-4-2}), we have that $d'\sim v_1$. Now, $\{d,d',v_1,d'',v\}$ induces a fork, a contradiction. So, $T'$ is complete to $V(C)$.   

Therefore, $V(C)\cup T$ is a homogeneous set as $V(C)\cup T$ is complete to $X\cup T'\cup Z''$ and anticomplete to $A$,  contradicts Lemma \ref{pb-2}. This proves Lemma~\ref{2-5}. \qed

\medskip

\emph{ Proof of Theorem~\ref{fork}.}  Suppose to its contrary. Let $G$ be a minimal nonperfectly divisible (fork, odd balloon)-free graph. It is certain that $G$ is connected.

By Lemma~\ref{pb-1}, $G[M(v)]$ is not perfect for any vertex $v$ of $G$. According to the strong perfect graph theorem and Lemma~\ref{2-3}, $G[M(v)]$ contains an odd hole, say $C=v_1v_2\dots v_nv_1$ such that $n\ge5$ and $n$ is odd. Let $A=M(C)$ and $u$ be a vertex in $N(A)$. Then $u$ has a neighbor in $V(C)$ and a neighbor in $A$. By Lemma~\ref{2-4}, we have that $N(u)\cap V(C)=\{v_j, v_{j+1}\}$, for some $j\in \{1,2,\dots,n\}$, if $N(u)\cap V(C)\ne V(C)$.

Let $W_1=\{v\in N(A)| N(v)\cap V(C)=V(C)\}$, let $W_2=N(A)\setminus W_1$. Since $G$ is odd balloon-free, we have that $W_2=\emptyset$, which contradicts Lemma~\ref{2-5}. This completes the proof of Theorem~\ref{fork}.\qed

\section{(fork, gem)-free graphs }

In this section, we consider (fork,gem)-free graphs, and prove Theorem~\ref{gem}.  A vertex of a graph is {\em bisimplicial} if its neighborhood is the union of two cliques. We need the following Lemma \ref{bisimplicial}, which establishes a connection between odd balloons and bisimplicial vertices in a (fork,gem)-free graph. Let $B$ be an induced odd balloon in $G$, and $u\in V(B)$. We call $u$ be the {\em center} of $B$ if $d_B(u)=3$ and $u$ has a neighbor of degree 1 in $B$, and call $B$ a {\em minimum induced odd balloon} if $|V(B)|$ is minimum among all induced odd balloons of $G$.

\begin{lemma}\label{bisimplicial}
	Let $G$ be a $($fork,gem$)$-free graph. Then $G$ has a bisimplicial vertex if $G$ contains an induced odd balloon.
\end{lemma}

\pf Let $B$ be a minimum induced odd balloon of $G$, and $u$ be the center of $B$. Let $C=v_1v_2\dots v_nv_1$ be the odd hole of $B$ such that $n\ge5$ and $n$ is odd, and $s$ be the vertex of degree 1 of $B$. Without loss of generality, we may assume that $u$ is complete to $\{v_1,v_2\}$ and anticomplete to $\{v_3,v_4,\dots,v_n\}$.

Since $G$ is gem-free, we have that

\begin{equation}\label{gem-1}
\mbox{no vertex of $V(G)\setminus V(C)$ may have four consecutive neighbors in $V(C)$.}
\end{equation}

Let $A$ denote the set of vertices which are complete to $\{v_1,v_2\}$ and anticomplete to $\{v_3,v_4\}$. If $A$ has two nonadjacent vertices $a_1$ and $a_2$, then $\{a_1,a_2,v_2,v_3,v_4\}$ induces a fork, a contradiction. So, we have that

\begin{equation}\label{gem-2}
\mbox{$A$ is a clique.}
\end{equation}

Let $D$ denote the set of vertices in $N(u)\setminus\{v_1,v_2\}$ which are anticomplete to $V(C)$. Clearly, $s\in D$.  If $D$ has two nonadjacent vertices $d_1$ and $d_2$, then $\{d_1,d_2,u,v_2,v_3\}$ induces a fork, a contradiction. Therefore,

\begin{equation}\label{gem-3}
\mbox{$D$ is a clique.}
\end{equation}

Let $U$ be the set of vertices in $N(u)$ which are complete to $D$. We will prove that $N(u)\setminus (U\cup D)$ and $U\cup D$ are both cliques in the following next two claims.

\begin{claim}\label{claim-1}
	Let $T$ be the set of verices in $N(u)\setminus\{v_1,v_2\}$ which are nonadjacent to $s$. Then $T$ is a clique.
\end{claim}

\pf Firstly, we prove that

\begin{equation}\label{gem-4}
\mbox{all vertices of $T$ are complete to $\{v_1,v_2\}$.}
\end{equation}

Suppose to its contrary. Let $t$ be a vertex of $T$ such that $t\not\sim v_2$ by symmetry. To forbid a fork $\{s,t,u,v_2,v_3\}$, we have that $t\sim v_3$. Suppose $t\sim v_n$, $\{s,u,t,v_3,v_n\}$ induces a fork, a contradiction. This proves $t\not\sim v_n$. To forbid a fork $\{s,t,u,v_1,v_n\}$, we have that $t\sim v_1$. If $t\sim v_{n-1}$, then $\{s,u,t,v_2,v_{n-1}\}$ induces a fork, a contradiction. This proves $t\not\sim v_{n-1}$. But now, $\{t,v_1,v_2,v_{n-1},v_n\}$ induces a fork, a contradiction. This proves (\ref{gem-4}).

Next, we prove that

\begin{equation}\label{gem-5}
\mbox{if $t\in T$, then $N(t)\cap V(C)\in \{\{v_1,v_2\},\{v_1,v_2,v_3\},\{v_1,v_2,v_n\}\}$}
\end{equation}

Suppose to its contrary. By (\ref{gem-4}), we have that $t$ is complete to $\{v_1,v_2\}$, which implies that there exists an integer $j\in\{4,\dots,n-1\}$ such that $t\sim v_j$. By (\ref{gem-1}), we have that $t\not\sim v_3$ or $t\not\sim v_n$. If $t$ has exactly one neighbor in $\{v_3,v_n\}$, we suppose by symmetry that $t\sim v_3$. By (\ref{gem-1}), we have that $t\not\sim v_4$. Let $j'\in\{5,\dots,n-1\}$ be the smallest integer such that $t\sim v_{j'}$. Then $\{t,v_3,v_{j'},u,s\}$ induces a fork, a contradiction. So, we suppose that  $t\not\sim v_3$ and $t\not\sim v_n$.

If $n=5$, then $t\not\sim v_3, t\not\sim v_5$ and $t\sim v_4$, which implies that $\{v_3,v_4,v_5,t,u\}$ induces a fork, a contradiction. Hence, we may assume that $n\ge7$.

If $t$ is anticomplete to $V(C)\setminus\{v_1,v_2,v_3,v_j,v_n\}$, then $\{v_{j-1},v_j,v_{j+1},t,u\}$ induces a fork, a contradiction. So, there exists an integer $i\in\{4,\dots,n-1\}\setminus\{j\}$ such that $t\sim v_i$. If $i\not\in \{j-1, j+1\}$, then $\{v_i,v_j,t,u,s\}$ induces a fork, a contradiction. This proves that $i\in \{j-1, j+1\}$. If $t$ has exactly one neighbor in $\{v_{j-1},v_{j+1}\}$, we suppose by symmetry that $t\sim v_{j-1}$, then $N(t)\cap \{v_4,\dots,v_{n-1}\}=\{v_{j-1}, v_j\}$, which implies that $\{v_j,\dots,v_n,v_1,t,u,s\}$ or $\{v_2,v_3,\dots,v_{j-1},t,u,s\}$ induces a smaller odd balloon as $n\ge7$, a contradiction. Therefore, $t\sim v_{j-1}$ and $t\sim v_{j+1}$, and thus $\{v_{j-1},v_{j+1},t,u,s\}$ induces a fork, a contradiction. This proves (\ref{gem-5}).

Let $T_1$ be the subset of $T$ such that $N(T_1)\cap V(C)=\{v_1,v_2\}$, $T_2$ be the subset of $T$ such that $N(T_2)\cap V(C)=\{v_1,v_2,v_3\}$, and let $T_3$ be the subset of $T$ such that $N(T_3)\cap V(C)=\{v_1,v_2,v_n\}$. Then, $T=T_1\cup T_2\cup T_3$ by (\ref{gem-5}). By (\ref{gem-2}), $T_1$ is a clique. If $T_2$ is not a clique, let $t',t''\in T_2$ be two nonadjacent vertices, then $\{t',t'',v_3,v_4,v_5\}$ induces a fork, a contradiction. Therefore, $T_2$ is a clique. Similarly, $T_3$ is a clique.

Suppose that $T$ is not a clique. Let $t_1$ and $t_2$ be two nonadjacent vertices of $T$. Since $T_1,T_2$ and $T_3$ are all cliques, we have that $t_1$ and $t_2$ belong to different elements of $\{T_1,T_2,T_3\}$. If $t_1\in T_2$ and $t_2\in T_3$, then $\{t_2,t_1,u,s,v_n\}$ induces a fork. Otherwise, we suppose by symmetry that $t_1\in T_1$ and $t_2\in T_2$, then $\{t_1,t_2,v_1,v_n,v_3\}$ induces a fork. Both are contradictions. This proves Claim~\ref{claim-1}.\qed

\begin{claim}\label{claim-2}
	$U$ is a clique.
\end{claim}

\pf Suppose to its contrary. Let $u_1$ and $u_2$ be two nonadjacent vertices of $U$. By the definition of $U$, $u_1$ and $u_2$ are both complete to $D$. We will prove that

\begin{equation}\label{gem-6}
\mbox{$\{u_1,u_2\}$ is anticomplete to $\{v_3,v_n\}$.}
\end{equation}

If it is not the case, we may assume  by symmetry that $u_2\sim v_3$. Since $u_2\sim s$ and $s$ is anticomplete to $V(C)$, by Lemma \ref{2-4} and (\ref{gem-1}), we have that $N(u_2)\cap V(C)\subseteq \{\{v_2,v_3\},\{v_3,v_4\}\}$.

So, $u_2\not\sim v_1$, which implies that $u_1\sim v_1$ or $u_1\sim v_3$ to forbid a fork on $\{v_1,u_1,u,u_2,v_3\}$. Suppose $u_1\sim v_1$. Then $u_1\sim v_2$ to forbid a gem on $\{v_1,v_2,u_1,u,s\}$, and consequently, $v_2\sim u_2$ to forbid a gem on $\{v_2,u_1,u,u_2,s\}$, which leads to a contradiction as $\{v_2,v_3,u_2,u,s\}$ induces a gem. So $u_1\not\sim v_1$, and thus $u_1\sim v_3$. Since $u_1\sim s$ and $s$ is anticomplete to $V(C)$, by Lemma \ref{2-4} and (\ref{gem-1}), we have that $N(u_1)\cap V(C)\subseteq \{\{v_2,v_3\},\{v_3,v_4\}\}$.

Therefore, $\{u_1,u_2\}$ is anticomplete to $\{v_1,v_n\}$. But now, $\{u_1,u_2,u,v_1,v_n\}$ induces a fork, a contradiction. This proves (\ref{gem-6}).

Next, we prove that

\begin{equation}\label{gem-7}
\mbox{$\{u_1,u_2\}$ is anticomplete to $\{v_1,v_2\}$.}
\end{equation}

Suppose to its contrary, we may assume that $u_2\sim v_2$ by symmetry. Since $u_2\sim s$ and $s$ is anticomplete to $V(C)$, by Lemma \ref{2-4} and (\ref{gem-1}), we have that $N(u_2)\cap V(C)\subseteq \{\{v_1,v_2\},\{v_2,v_3\}\}$. Since $u_2\not\sim v_3$ by (\ref{gem-6}), we have that $N(u_2)\cap V(C)= \{v_1,v_2\}$. Moreover,
$u_1\sim v_2$ to forbid a gem on $\{v_2,u_1,u,u_2,s\}$. Similarly, we have that $N(u_1)\cap V(C)=N(u_2)\cap V(C)= \{v_1,v_2\}$. Now, $\{u_1,u_2,v_2,v_3,v_4\}$ induces a fork, a contradiction. This proves (\ref{gem-7}).

By (\ref{gem-6}) and (\ref{gem-7}), $\{u_1,u_2\}$ is anticomplete to $\{v_1,v_2, v_3,v_n\}$, which implies that $\{u_1,u_2,u,v_2,v_3\}$ induces a fork, a contradiction. This proves Claim~\ref{claim-2}.\qed

\medskip

Let $T^+$ be the set of vertices in $N(u)\setminus (D\cup \{v_1,v_2\})$ which are not complete to $D$. We will prove that

\begin{equation}\label{gem-8}
\mbox{$T^+$ is a clique.}
\end{equation}

Suppose that $T^+$ is not a clique. Let $t_3$ and $t_4$ be two nonadjacent vertices of $T^+$. If $D$ has a vertex, say $d$, that is nonadjacent to both $t_3$ and $t_4$, by replacing $s$ with $d$ in Claim~\ref{claim-1}, we have that $t_3\sim t_4$, a contradiction. Therefore, there exist $d_3,d_4\in D$ such that $t_3\sim d_3, t_3\not\sim d_4, t_4\sim d_4$ and $t_4\not\sim d_3$. Since $D$ is a clique by (\ref{gem-3}), we have that $d_3\sim d_4$, which implies that $\{u,t_3,t_4,d_3,d_4\}$ induces a gem, a contradiction. This proves (\ref{gem-8}).

By (\ref{gem-4}), $T$ is complete to $\{v_1,v_2\}$. Since each vertex $t^+$ of $T^+$ has a nonneighbor in $D$, say $d^+$, by replacing $s$ with $d^+$ in Claim~\ref{claim-1}, we have that $t^+$ is complete to $\{v_1,v_2\}$, and thus $T^+$ is complete to $\{v_1,v_2\}$. Therefore, by (\ref{gem-8}), $N(u)\setminus (U\cup D)=T^+\cup\{v_1,v_2\}$ is a clique. By the definition of $U$ and Claim~\ref{claim-2}, $U\cup D$ is a clique as $D$ is a clique. That is to say, $u$ is a bisimplicial vertex of $G$. This completes the proof of Lemma~\ref{bisimplicial}.
\qed

\medskip

\emph{ Proof of Theorem~\ref{gem}.} If $G$ has an induced odd balloon, then $G$ has a bisimplicial vertex, say $u$ such that $d(u)\le 2\omega(G)-2$ by Lemma~\ref{bisimplicial}. By induction, we may suppose that $\chi(G-u)\le\frac{1}{2}\omega^2(G-u)+\frac{1}{2}\omega(G-u)\le\frac{1}{2}\omega^2(G)+\frac{1}{2}\omega(G)$. Since for any $\omega(G)\ge1, \frac{1}{2}\omega^2(G)+\frac{1}{2}\omega(G)-(2\omega(G)-2)\ge1$, we can take any $\chi(G-u)$-coloring of $G-u$ and extend it to a $(\frac{1}{2}\omega^2(G)+\frac{1}{2}\omega(G))$-coloring of $G$.

If $G$ has no induced odd balloon, then $G$ is perfectly divisible by Theorem~\ref{fork}, and hence $\chi(G)\le \frac{1}{2}\omega^2(G)+\frac{1}{2}\omega(G)$. This completes the proof of Theorem~\ref{gem}.
\qed

\section{(fork, butterfly)-free graphs }

In this section, we consider (fork,butterfly)-free graphs, and prove Theorem~\ref{butterfly}.  A vertex of a graph is {\em trisimplicial} if its neighborhood is the union of three cliques. With a similar technique used in the proof of Theorem~\ref{gem}, we can establish a connection between odd balloons and trisimplicial vertices in a (fork,butterfly)-free graph. It is known that for a fork-free graph $G$, $\chi(G)\le3$ if $\omega(G)=2$, and  $\chi(G)\le4$ if $\omega(G)=3$\cite{RB93}. Therefore, we may suppose that $\omega(G)\ge4$ in this section.

\begin{lemma}\label{trisimplicial}
	Let $G$ be a $($fork,butterfly$)$-free graph. Then $G$ has a trisimplicial vertex if $G$ contains an induced odd balloon.
\end{lemma}
\pf Let $B$ be a minimum induced odd balloon of $G$, and $u$ be the center of $B$. Let $C=v_1v_2\dots v_nv_1$ be the odd hole of $B$ such that $n\ge5$ and $n$ is odd, and $s$ be the vertex of degree 1 of $B$. Without loss of generality, we may assume that $u$ is complete to $\{v_1,v_2\}$ and anticomplete to $\{v_3,v_4,\dots,v_n\}$.

Let $A$ denote the set of vertices which are complete to $\{v_1,v_2\}$ and anticomplete to $\{v_3,v_4\}$. With the same argument as used to prove (\ref{gem-2}), we can show that

\begin{equation}\label{butterfly-0}
\mbox{$A$ is a clique.}
\end{equation}

Let $D$ denote the set of vertices in $N(u)\setminus\{v_1,v_2\}$ which are anticomplete to $V(C)$. It is clear that $s\in D$. With the same argument as used to prove (\ref{gem-3}), we can show that $D$ is a clique. In particular, if $D$ has two distinct vertices, say $d_1$ and $d_2$, then $\{v_1,v_2,u,d_1,d_2\}$ induces a butterfly, a contradiction. This proves that $D=\{s\}$.

Let $T$ be the set of vertices in $N(u)\setminus\{v_1,v_2\}$ which are nonadjacent to $s$. With the same argument as used to prove (\ref{gem-4}), we can show that

\begin{equation}\label{butterfly-1}
\mbox{all vertices of $T$ are complete to $\{v_1,v_2\}$.}
\end{equation}

Next, we will prove that

\begin{equation}\label{butterfly-2}
\mbox{no vertex of $T$ may have four consecutive neighbors in $V(C)$.}
\end{equation}

If it is not true, let $t$ be a vertex in $T$ that has four consecutive neighbors $\{v_i,v_{i+1},v_{i+2},v_{i+3}\}$ in $V(C)$. By (\ref{butterfly-1}), we have that $t$ is complete to $\{v_1,v_2\}$. By the choice of $t$, we may suppose that $|\{v_i,v_{i+1},v_{i+2},v_{i+3}\}\cap \{v_1,v_2\}|\in\{0,2\}$.

If $\{v_i,v_{i+1},v_{i+2},v_{i+3}\}\cap \{v_1,v_2\}=\emptyset$, we may suppose that $i\ge4$ and $i+3\le n-1$, then $\{v_1,v_{i+1},v_{i+2},t,u\}$ induces a butterfly, a contradiction. So, $|\{v_i,v_{i+1},v_{i+2},v_{i+3}\}\cap \{v_1,v_2\}|=2$, which implies that $i\in\{1,n-1,n\}$. If $i=n$, then $\{v_3,v_n,t,u,s\}$ induces a fork as $t\not\sim s$ and $u$ is anticomplete to $\{v_3,v_n\}$, a contradiction. So, $i=1$ or $i=n-1$, we may assume by symmetry that $i=1$. Then $\{v_1,v_3,v_4,t,u\}$ induces a butterfly, a contradiction. Therefore $|\{v_i,v_{i+1},v_{i+2},v_{i+3}\}\cap \{v_1,v_2\}|\ne2$. This proves (\ref{butterfly-2}).

By replacing (\ref{gem-2}) and (\ref{gem-1}) with (\ref{butterfly-0}) and (\ref{butterfly-2}), respectively, we can prove (\ref{butterfly-3}) and (\ref{butterfly-4}), with the same arguments as that used in proving  Claim~\ref{claim-1}.

\begin{equation}\label{butterfly-3}
\mbox{If $t\in T$, then $N(t)\cap V(C)\in \{\{v_1,v_2\},\{v_1,v_2,v_3\},\{v_1,v_2,v_n\}\}$}
\end{equation}

and

\begin{equation}\label{butterfly-4}
\mbox{$T$ is a clique.}
\end{equation}

Let $U$ be the set of vertices in $N(u)$ which are adjacent to $s$. If $u'$ is a vertex in $U$ which is complete to $V(C)$, then $\{v_3,v_4,u',u,s\}$ induces a butterfly. If $u''$ is a vertex in $U$ which is anticomplete to $\{v_1,v_2\}$, then $\{v_1,v_2,u'',u,s\}$ induces a butterfly. Both are contradictions. So, we have that

\begin{equation}\label{butterfly-5}
\mbox{no vertex of $U$ is complete to $V(C)$, and no vertex of $U$ is anticomplete to $\{v_1,v_2\}$.}
\end{equation}

Let $U_1$ be the set of vertices in $U$ which are complete to $\{v_1,v_2\}$, and let $U_2=U\setminus U_1$. Since each vertex of $U_1$ is adjacent to $s$, by Lemma \ref{2-4} and (\ref{butterfly-5}), $N(U_1)\cap V(C)=\{v_1,v_2\}$. If $U_1$ has two nonadjacent vertices $u_1'$ and $u_1''$, then $\{u_1',u_1'',v_2,v_3,v_4\}$ induces a fork, a contradiction, which implies that $U_1$ is a clique. Next, we will prove that

\begin{equation}\label{butterfly-6}
\mbox{$T\cup U_1$ is a clique.}
\end{equation}

Suppose to its contrary that there exist two nonadjacent vertices $t$ and $u_1$ in $T\cup U_1$. Since $T$ and $U_1$ are both cliques, we may assume that $t\in T$ and $u_1\in U_1$. By (\ref{butterfly-3}), $N(t)\cap V(C)\subseteq \{\{v_1,v_2\},\{v_1,v_2,v_3\},\{v_1,v_2,v_n\}\}$. Moreover, $N(u_1)\cap V(C)=\{v_1,v_2\}$ as $N(U_1)\cap V(C)=\{v_1,v_2\}$ by Lemma \ref{2-4} and (\ref{butterfly-5}). If $N(t)\cap V(C)\subseteq \{\{v_1,v_2\},\{v_1,v_2,v_n\}\}$, then $\{t,u_1,v_2,v_3,v_4\}$ induces a fork. If $N(t)\cap V(C)=\{v_1,v_2,v_3\}$, then $\{t,u_1,v_1,v_{n-1},v_n\}$ induces a fork. Both are contradictions. This proves (\ref{butterfly-6}).

For $i\in\{1,2\}$, let $U_2^i=U_2\cap N(v_i)$. Since no vertex of $U$ is anticomplete to $\{v_1,v_2\}$ by (\ref{butterfly-5}), we have that $v_2$ is anticomplete to $U_2^1$ and $v_1$ is anticomplete to $U_2^2$, and $U_2=U_2^1\cup U_2^2$. Now, we will prove that

\begin{equation}\label{butterfly-7}
\mbox{$U_2^1$ and $U_2^2$ are both cliques.}
\end{equation}

If it is not the case, we may assume by symmetry that $U_2^1$ is not a clique, then there exist two nonadjacent vertices $w$ and $w'$ in $U_2^1$. Since $w$ and $w'$ are both adjacent to $v_1$ and $s$, by Lemma \ref{2-4} and (\ref{butterfly-5}), we have that $N(w)\cap V(C)\in\{\{v_1,v_2\},\{v_1,v_n\}\}$ and $N(w')\cap V(C)\in\{\{v_1,v_2\},\{v_1,v_n\}\}$. Therefore, $N(w)\cap V(C)=N(w')\cap V(C)=\{v_1,v_n\}$ as $v_2$ is anticomplete to $U_2^1$. But now, $\{w,w',v_1,v_2,v_3\}$ induces a fork, a contradiction. This proves (\ref{butterfly-7}).

By the definition of $U_1$, $U_1$ is complete to $\{v_1,v_2\}$. By (\ref{butterfly-1}), $T$ is also complete to $\{v_1,v_2\}$. So, $T\cup U_1\cup\{v_1,v_2\}$ is a clique as $T\cup U_1$ is a clique by (\ref{butterfly-6}). By the definition of $U$, we have that $s$ is complete to $U_2^1$, which implies that $U_2^1\cup\{s\}$ is a clique as $U_2^1$ is a clique by (\ref{butterfly-7}). Since $U_2\cup\{s\}=(U_2^1\cup\{s\})\cup U_2^2=N(u)\setminus (T\cup U_1\cup\{v_1,v_2\})$, we have that $u$ is a trisimplicial vertex of $G$ as $U_2^1\cup\{s\},U_2^2$ and $T\cup U_1\cup\{v_1,v_2\}$ are all cliques. This completes the proof of Lemma~\ref{trisimplicial}.
\qed

\medskip

\emph{ Proof of Theorem~\ref{butterfly}.} If $G$ has an induced odd balloon, then $G$ has a trisimplicial vertex, say $u$ such that $d(u)\le 3\omega(G)-3$ by Lemma~\ref{trisimplicial}. By induction, we may suppose that $\chi(G-u)\le\frac{1}{2}\omega^2(G-u)+\frac{1}{2}\omega(G-u)\le\frac{1}{2}\omega^2(G)+\frac{1}{2}\omega(G)$. Since for any $\omega(G)\ge4, \frac{1}{2}\omega^2(G)+\frac{1}{2}\omega(G)-(3\omega(G)-3)\ge1$, we can take any $\chi(G-u)$-coloring of $G-u$ and extend it to a $(\frac{1}{2}\omega^2(G)+\frac{1}{2}\omega(G))$-coloring of $G$.

If $G$ has no induced odd balloon, then $G$ is perfectly divisible by Theorem~\ref{fork}, and hence $\chi(G)\le \frac{1}{2}\omega^2(G)+\frac{1}{2}\omega(G)$. This completes the proof of Theorem~\ref{butterfly}. \qed


\begin{thebibliography}{9999}
	
	\bibitem{BM08} J. A. Bondy and U. S. R. Murty. Graph Theory, Springer, New York, 2008.
	
	
	
	\bibitem{CCS20} M. Chudnovsky, L. Cook and P. Seymour. Excluding the fork and antifork. Discrete Mathematics  343 (2020), 111786.
	
	\bibitem{CHKK21}M. Chudnovsky, S. Huang, T. Karthick and J. Kaufmann. Square-free graphs with no induced fork. The Electronic Journal of Combinatorics 28(2) (2021) \#P2.20.
	
	\bibitem{CRSR06} M. Chudnovsky, N. Robertson, P. Seymour and R. Thomas. The strong perfect graph theorem. Annals of mathematics 164(1) (2006) 51-229.
	
	\bibitem{CS19} M. Chudnovsky and V. Sivaraman. Perfect divisibility and 2-divisibility. Journal of Graph Theory 90 (2019) 54-60.
	
	
	
	\bibitem{gyarfas1} A. Gy\'{a}rf\'{a}s. On Ramsey covering-numbers. {\it Colloquia Mathematic Societatis J\'{a}nos Bolyai 10, Infinite and Finite Sets.} North-Holland/American Elsevier, New York. (1975) 801-816.
	
	\bibitem{HCT18} C.~T. Ho\'{a}ng. On the structure of (banner, odd hole)-free graphs. Journal of Graph Theory 89 (2018) 395--412.
	
	\bibitem{KKS22} T. Karthick, J. Kaufmann and V. Sivaraman. Coloring graph classes with no induced fork via perfect divisibility. The Electronic Journal of Combinatorics 29(3) (2022) \#P3.19.
	
	\bibitem{KP94} H.~A. Kierstead and S.~G. Penrice. Radius two trees specify $\chi$-bounded classes. Journal of Graph Theory 18(2) (1994) 119-129.
	
	\bibitem{LSWY21} X. Liu, J. Schroeder, Z. Wang and X. Yu. Polynomial $\chi $-binding functions for $ t $-broom-free graphs. arXiv preprint arXiv:2106.08871 (2021).
	
	\bibitem{RB93} B. Randerath. The Vizing bound for the chromatic number based on forbidden pairs. PhD thesis, RWTH Aachen, Shaker Verlag, 1993.
	
	\bibitem{RS04} B. Randerath and I. Schiermeyer. Vertex colouring and forbidden subgraphs - a survey. Graphs and Combinatorics 20(1) (2004) 1-40.
	
	
	
	
\end{thebibliography}
\end{document}